\documentstyle[12pt]{article} 
\begin{document}
\setcounter{page}{1} 
\begin{center}
{\Large\bf {\bf The M\"obius function on the lattice of normal subgroups}} 
\\by 
\\\ Gregory M.  
Constantine {\Large\ } 
\end{center}
\medskip
\centerline{Department of Mathematics}

\centerline{The University of Pittsburgh}

\centerline{Pittsburgh, PA 15260}
\vskip1cm
\centerline{\bf ABSTRACT}
 
\medbreak\noindent By studying lattices of normal 
subgroups, especially those of the socle and radical, an 
expression is obtained for the minimal number of 
conjugacy classes required to generate a group.  This 
number is shown to be captured by the character table.  
The M\"obius function is then used to extract information 
on the faithful irreducible representations of a group.  

\medbreak\noindent {\em Key words and phrases:\/}  Simple 
module, character, faithful representation, socle, radical.  
\medbreak\noindent {\em Proposed running head:\/}  M\"obius 
function on normal subgroups 
\vskip2cm 
\footnoterule\noindent Funded under NIH grants 
P50-GM-53789, RO1-HL-076157 and an IBM shared 
University Research Award 

\noindent gmc@math.pitt.edu 
  
\newpage
\noindent

\section*{\Large\bf1.  The lattice of normal subgroups} 
\hspace{\parindent} All groups considered in this paper 
are finite.  For a subgroup $H$ of $G$ the set of subgroups 
of $H$ that are normal in $G$, with inclusion on subgroups 
as partial order, forms a lattice $L_G(H)$ with meet $T\cap K$ 
and join $TK$.  We write $L(G)$ for $L_G(G).$ In accordance 
with the theory of linear representations over the 
complex field, it helps to visualize this lattice as the 
intersection of kernels of irreducible representations of 
$G,$ especially in the latter part of the paper; cf.  
Alperin and Bell (1995, p 150).  Silcock (1977) proved that 
every finite distributive lattice is isomorphic to the 
lattice of normal subgroups of a group.  The group 
involved is a wreath product of a family of non-abelian 
simple groups indexed by a partially ordered set.  
Making use of well-known results of Weisner (1935), Hall 
(1933, 1935) and Rota (1964), and the structure of the 
minimal normal subgroups viewed as irreducible 
$G-$modules over a prime-order field, we first obtain an 
explicit expression for the M\"obius function of $L(G)$.  
Basics on M\"obius functions on partially ordered sets are 
found in Constantine (1987, Chapter 9).  
\vskip.2cm 
\noindent\ In the context of this paper, Hall's result 
states that $\mu (H,G)=0$ unless $H$ is the intersection of 
maximal normal subgroups of $G$, which are the coatoms 
of $L(G)$.  Dually, it informs us that $\mu (1,K)=0$ unless $K$ is 
the join of minimal normal subgroups of $G,$which are the 
atoms of $L(G).$ We shall derive a general expression for 
$\mu (K,H)$ with $K$ and $H$ normal subgroups of $G;$ it is, 
however, useful to first examine certain special 
segments of the lattice $L(G),$ since they are directly 
involved in many inversion formulae.  As indicated, the 
segments $[R,G]$ and $[1,S]$ are particularly interesting, 
with $R$ signifying the intersection of all maximal normal 
subgroups of $G,$ and $S$ being the subgroup generated by 
all minimal normal subgroups of $G.$ The subgroup $R$ is 
called the {\em radical\/} and $S$ is the {\em socle\/} of $G.$ \medbreak It 
is well-known that any minimal normal subgroup is a 
direct product of isomorphic simple groups, cf.  Alperin 
and Bell (1995, p 94).  It follows that $\mu (1,G)=0,$ unless 
$G$ is a direct product of simple groups.  
\vskip.3cm 
An explicit form of $\mu (1,G)$ when $G$ is a direct product of 
simple groups is obtained by using a basic property of 
the lattice of normal subgroups of such a group, which 
we now describe.  
\vskip.2cm 
For $(X,\leq )$ and $(Y,\leq )$ partially ordered sets (posets, for 
short), we define a poset $(X\times Y,\leq )$ by setting 
$(x_1,y_1)\leq (x_2,y_2)$ if $x_1\leq x_2$ and $y_1\leq y_2.$ The newly defined 
poset is said to be the {\em product\/} of the two initial posets.  
It is well-known, see Constantine (1987, p 383), that the 
M\"obius function of a product of posets is the pointwise 
product of the respective M\"obius functions; specifically, 
with obvious notation, 
$\mu ((x_1,y_1),(x_2,y_2))=\mu_X(x_1,x_2)\mu_Y(y_1,y_2).$ Seminal work on 
group actions on posets is found in Stanley (1982).  
Coming back to the lattice of normal subgroups, in 
general we cannot say that $L(G\times H)=L(G)\times L(H).$ Observe 
that in order to have $L(G\times H)=L(G)\times L(H)$ we should 
ensure that every normal subgroup of $G\times H$ is the direct 
product of a normal subgroup of $G$ and a normal 
subgroup of $H.$ Miller (1975) gives a characterization of 
direct products of groups with this property.  
Specifically, he shows that the lattice of normal 
subgroups of a group $G$ is a non-trivial direct product if 
and only if $G=G_1\times G_2$ with $G_1\neq 1\neq G_2$ and either for at 
least one factor $G_i[G_i,N_i]=N_i$ for every normal 
subgroup $N_i$ of $G_i$, or the elements of $N_1/[G_1,N_1]$ and 
$N_2/[G_2,N_2]$ have finite relatively prime orders for 
every normal subgroup $N_1$ of $G_1$ and $N_2$ of $G_2$.  Our 
interests in this paper are focussed on special kinds of 
direct products, associated with the socle of a group $G$, 
and with quotients to the radical of $G$.  
\vskip.2cm 
\noindent {\bf Lemma 1 }
\vskip.1cm 
\noindent(a) {\em If} $G_1$ {\em and} $G_2$ {\em are groups of coprime orders, }
{\em then} $L(G_1\times G_2)=L(G_1)\times L(G_2).$ 
\vskip.1cm 
\noindent(b) {\em If S is non-Abelian simple and G is }
{\em arbitrary, then} $L(S\times G)=L(S)\times L(G).$ 
\vskip.2cm 
\noindent {\bf Proof:}  (a) Let $N\unlhd (G_1\times G_2).$ We need to show 
that $N=N_1\times N_2$, for some $N_i\unlhd G_i.$ In fact we can select 
$N_i=N\cap G_i.$ Indeed, the $N_i$ so defined is normal in $G_i$.  
Clearly $N_1\times N_2\leq N.$ Take $n\in N,$ $n=(g_1,g_2)$, with $
g_i\in G_i.$ 
We are done if we prove that $g_i\in N_i.$ Since $g_1$ and $g_2$ 
are elements of coprime orders, we have 
$<n>=<(g_1,g_2)>=<g_1>\times <g_2>.$ It follows that 
$g_i\in <n>\leq N,$ and therefore that $g_i\in N\cap G_i=N_i.$ 
\vskip.1cm 
(b) For $H\unlhd (S\times G)$ define $H_S=\{s\in S:$ there exists $
g\in G$ 
such that $(s,g)\in H$$\}.$ Similarly, $H_G=\{g\in G:$ there exists 
$s\in S$ such that $(s,g)\in H$$\}.$ We have $H_S\unlhd S$ and $H_
G\unlhd G.$ By 
the simplicity of $S,$ we know that $H_S$ is either 1 or $S.$ 
Consider $H\cap S;$ it is either 1 or $S.$ If it is 1, then 
$[S,H]\leq S\cap H=1,$ since both $S$ and $H$ are normal in 
$S\times G.$ Consequently, for all $(x,y)\in H$ and all $(s,1)\in 
S$ we 
have $(x,y)=(x,y)^{(s,1)}=(x^s,y).$ It follows that $H_S$ is in 
the center of $S$ and is therefore 1.  This implies $H\unlhd G,$ 
and we conclude that $H=1\times H.$ Assume now that 
$H\cap S=S.$ In this case $S\leq H,$ and therefore $H=S\times H_G
.$ 
This ends the proof.  
\vskip.2cm 
We can now give an explicit form for $\mu (1,G)$ when $G$ is 
a direct product of simple groups.  
\vskip.2cm 
\noindent {\bf Lemma 2 }
\vskip.1cm 
\noindent(a) {\em If G is not a direct product of simple }
{\em groups, then} $\mu (1,G)=0.$ 
\vskip.1cm 
\noindent(b) {\em If} $G=S_1^{s_1}\times\cdots\times S_m^{s_m}\times 
C_{p_1}^{a_1}\times\cdots\times C_{p_n}^{a_n}${\em , with} $S_i$ 
{\em nonisomorphic non-Abelian simple groups and} $C_{p_j}$ 
{\em nonisomorpic cyclic groups of prime order, then }

\[\mu (1,G)=(-1)^{\sum_{i=1}^ms_i}\prod_{j=1}^n(-1)^{a_j}p_j^{{{a_
j}\choose 2}}.\]
\noindent
{\bf Proof:} Part (a) had already been argued. As for part (b), 
observe that Lemma 2 is a
consequence of Lemma 1, since the lattice $L(G)$ is the 
product of the lattices written below 
\[L(G)=L(\prod_iS_i^{s_i}\times\prod_jC_{p_j}^{a_j})=\prod_iL(S_i
)^{s_i}\prod_jL(C_{p_j}^{a_j}).\]
\noindent
For a simple group $S$ we have 
$\mu (1,S)=-1,$ and it is well-known that $\mu (1,C_p^a)=(-1)^ap^{{
a\choose 2}}.$
\noindent
[We note in passing that clearly $L(C_p^a)\neq L(C_p)^a.$] 
\vskip.3cm
Denote by $R(G)$ the radical of $G.$ The subgroup $R(G),$
written simply as $R$ when the context permits, is a
characteristic subgroup of $G$ and has properties similar 
to the Frattini subgroup. Specifically as follows:
\vskip.2cm
\noindent
{\bf Lemma 3 }
\vskip.1cm
\noindent
(a) {\em Let Y be a set of conjugacy classes of G. Then, }
$<Y,R>=G$ {\em if and only if} $<Y>=G.$
\vskip.1cm
\noindent
(b) {\em For groups A and B, we have} $R(A\times B)=R(A)\times R(
B).$
\vskip.1cm
\noindent
(c) {\em Group G is a direct product of simple groups if }
{\em and only if} $R=1.$
\vskip.1cm
\noindent
(d) {\em The factor group G/R is a direct product of simple }
{\em groups, and R}
{\em is the smallest normal subgroup of G with this property.   }
\vskip.2cm
\noindent
{\bf Proof:} (a) Since conjugation is an automorphism, it is 
evident that a subgroup of $G$ generated 
by conjugacy classes is always normal in $G.$ If 
$<Y,R>=G$ but $<Y><G,$ then $<Y>\leq M,$ for $M$ some 
maximal normal subgroup of $G.$ Since also $R\leq M,$ we conclude
that $G=<Y,R>\leq M,$ which is a contradiction.
\vskip.1cm
(b) A maximal normal subgroup of $A\times B$ 
intersects precisely one of $A$ or $B$ in a maximal normal 
subgroup.
\vskip.1cm
(c) Denote by $M_1,M_2,\ldots$ the maximal normal subgroups 
of $G,$ and assume $R(G)=1.$ Consider $1\neq N\unlhd G.$ By
assumption there exists a maximal normal subgroup $M$ of 
$G$ which does not contain $N.$ Then 
\[N/(M\cap N)\cong NM/M=G/M,\]
\noindent
is a simple group, allowing us to conclude that $M\cap N$ is 
maximal normal in $N.$ It follows that 
$R(N)\leq\cap_i(M_i\cap N)=(\cap_iM_i)\cap N=1.$ Select $N=M_
1,$ and 
proceed by induction on the order of the group; we can 
therefore 
assume that $M_1$ is a direct product of simple groups, 
since $M_1\lhd G$ and $R(M_1)=1.$ To find a direct complement
to $M_1$ intersect just enough $M_is$, say $M_i$ with $i\in 
I$ (and 
note that index 1 is not in $I),$  
such that 
$\cap_{i\in I}M_i\neq 1$ but $M_1\cap (\cap_{i\in I}M_i)=1.$ Let $
H=\cap_{i\in I}M_i.$ We 
have $G=M_1\times H.$ Since $H$ is a proper normal subgroup 
of $G$ by induction $H$ is also a direct product of simple 
groups. 
\vskip.1cm
The converse is true since $R(S)=1$ for any simple group,
and by part (b) $R(\prod S_j)=\prod R(S_j)=1,$ with $S_j$ simple.
\vskip.1cm
(d) By the correspondence theorem we conclude that 
$R(G/R(G))=1.$ Part (c) assures us now that $G/R(G)$ is as
stated. On the other hand, if $G/H$ is a direct product of 
simple groups, then (c) tells us that $R(G/H)=1.$ Invoking
the
correspondence theorem again, we conclude that $H$ is the 
intersection of some maximal normal subgroups of $G;$ 
thus $R(G)\leq H.$ This ends the proof.
\vskip.3cm
In case $R\leq H$ we obtain
$\mu (H,G)=\mu (1,G/H),$ since the lattice of normal subgroups 
of $G$ that contain $H$ is isomorphic to the lattice of 
normal subgroups of $G/H$ by the correspondence 
theorem. Being a homomorphic image of $G/R,$ the 
group $G/H$ is in this case a product of 
simple groups; the value $\mu (H,G)$ is therefore 
known by Lemma 2(b). We summarize:
\vskip.2cm
\noindent
{\bf Lemma 4} {\em The M\"obius function} $\mu (H,G)=0,$ {\em unless} $
R\leq H,$ 
{\em in which case} $G/H$ {\em is a direct product of simple }
{\em groups and the value} $\mu (H,G)=\mu (1,G/H)$ {\em is found in }
{\em Lemma 2(b).}
\vskip.2cm
We stress that if 
$\Sigma (G)$ denotes the lattice of {\em all\/} subgroups of $G$, then we 
have the well beahaved property $\Sigma (G)\cap H=\Sigma (H).$ 
The lattice $L(G),$ however, does not possess 
this property. Indeed, for $H\unlhd G,$ in general 
$L(G)\cap H\neq L(H),$ since a normal subgroup of $H$ is not 
necessarily also normal in $G.$ Nevertheless we do have, 
again by the correspondence theorem, that  
$\mu (K,H)$ in the lattice $L(G)$ is equal to $\mu (1,H/K)$ in the 
lattice $L(G/K).$ It suffices therefore to find an explicit 
expression for $\mu (1,H)$ for $H\in L(G).$ We are thus led to 
study the lattice $L_G(S),$ where $S$ is the socle of $G.$ 
\vskip.5cm
View an elementary Abelian minimal normal
subgroup $M$ of $G$ as a simple $G-$module with conjugation 
as $G-$module action. Call two such minimal normal subgroups
$G-${\em isomorphic\/} if there exists a group 
isomorphism $f$ between them that preserves the 
$G-$module action; that is, $f(x^g)=f(x)^g,$ for all $g\in G.$ 
Partition the set of elementary Abelian 
minimal normal subgroups of $G$ into 
$G-$isomorphism classes. Denote these $G-$classes by 
${\bf A}_1,\ldots ,{\bf A}_a$.  Let $A_i$ be representatives from the
$G-$classes ${\bf A}_i$. The group generated by 
the elements of ${\bf A}_i$ is written as $S_i$. 
It follows that $S_i\cong A_i^{d_i},$ a
direct product of $d_i$ copies of $A_i.$ Let $S_{a+1},\ldots ,S_{
a+b}$ be the 
set of non-Abelian minimal normal subgroups of $G.$ Each 
$S_i,$ $a+1\leq i\leq a+b$, is a direct product of isomorphic non-Abelian simple 
groups. We can therefore express the socle $S$ of $G$ as 
the following direct product:\hfill\break
\centerline{$S=\prod_{i=1}^{a+b}S_i$.}
\vskip.2cm
Let $N\in L_G(S).$ Then $N=\prod_{i=1}^{a+b}(N\cap S_i).$ 
Indeed, we are done if $N=1;$ else $N$ contains a minimal 
normal subgroup $M$ of $G.$ Write $S=M\times C,$ and define 
$C_N=\{c\in C:(m,c)\in N,$ for some $m\in M$$\}$. It follows that 
$C_N$ is a normal subgroup of $G,$ and $N=M\times C_N.$ Since 
$C_N$ has lesser cardinality than $N$ we are done by 
induction on $\vert N\vert .$ To summarize:
\vskip.2cm
\noindent
{\bf Lemma 5} {\em The lattice} $L_G(S)$ {\em is equal to the direct }
{\em product} $\prod_{i=1}^aL_G(S_i)\prod_{j=a+1}^{a+b}L_G(S_j).$
\vskip.3cm
The lattice $L_G(S_j),$ $a+1\leq j\leq a+b,$ consists of just two elements $
\{1,S_j\}.$ 
By contrast, the structure of $L_G(S_i)$, $1\leq i\leq a,$ is more intricate. 
For a $G-$module $A,$ denote by $Hom_G(A,A)$ the set of 
$G-$homomorphisms from 
$A$ to $A$. We have the following result:
\vskip.2cm
\noindent
{\bf Lemma 6} {\em The lattice} $L_G(S_i)$ {\em with} $S_i\cong A_
i^{d_i}$ {\em is }
{\em isomorphic to the lattice of vector subspaces of a }
{\em vector space of dimension} $d_i$ {\em over the field }
$Hom_G(A_i,A_i).$
\vskip.2cm
\noindent
{\bf Proof:} By a result which Gorenstein (1968, page 79) attributes to 
J A Greene, without exact reference, it follows that there
are exactly $(q^{d_i}_i-1)(q_i-1)^{-1}$ minimal normal subgroups of $
G$ 
that are $G-$isomorphic to $A_i$, 
where $q_i$ denotes the number of $G-$module homomorphisms 
from $A_i$ to itself. By the same body of results, the 
number $q_i$ is in fact a power of the same prime that 
divides $\vert A_i\vert ,$ since the $G-$module homomorphisms from $
A_i$ 
to itself form in fact a finite field under the natural 
algebraic operations with such homomorphisms. We can 
now identify the normal subgroup $S_i\cong A_i^{d_i}$ of $G$ with a
vector space $V_i$ of dimension $d_i$ over the field $F_{q_i}$, and
establish a natural bijection between the lattice of 
normal subgroups of $G$ contained in $A_i^{d_i}$ and the lattice
of subspaces of $V_i$. Perhaps the simplest way to see the 
bijection is to first identify the minimal normal 
subgroups of $G$ isomorphic to $A_i$ with the 
one-dimensional subspaces of $V_i$. Subsequently view a 
normal subgroup $N$ of $G$ contained in $A_i^{d_i}$ as
generated by all the subgroups isomorphic to $A_i$ and normal 
in $G$ that it contains; the subspace of $V_i$ corresponding 
to $N$ is the subspace generated by all the corresponding 
one-dimensional subspaces to the subgroups isomorphic 
to $A_i$ contained in $N$ and normal in $G.$ Since the bijection 
is order preserving, it is a lattice isomorphism.
\vskip.3cm
The previous two lemmas yield the following information 
on the socle of a group:
\vskip.2cm
\noindent
{\bf Theorem 1} {\em Let the socle S of G be isomorphic to the }
{\em direct product} $\prod_{i=1}^aA_i^{d_i}\prod_{j=a+1}^{a+b}S_
j${\em , where} $\{A_i^{d_i}:1\leq i\leq a\}$ {\em are }
{\em representatives of the G-isomorphism classes of G with }
$A_i$ {\em elementary Abelian, and} $\{S_j:a+1\leq j\leq a+b\}$ {\em is the set of the }
{\em non-Abelian minimal normal subgroups of G. }
{\em Then} $L_G(S)\cong\prod_{i=1}^aL(V_i)\times \{0,1\}^b,$ {\em where} $
L(V_i)$ {\em is the lattice }
{\em of subspaces of a vector space of dimension} $d_i$ {\em over the }
{\em field} $F_{q_i},$ {\em with} $q_i=\vert Hom_G(A_i,A_i)\vert 
,$ {\em and} $\{0,1\}$ {\em denotes the }
{\em lattice with two elements.}
\vskip.3cm
Theorem 1 allows us to count the number of normal 
subgroups of $G$ of a given ''type'' contained in the socle 
of $G.$ It involves a product of Gaussian polynomials and 
binomial coefficients. This will be used in the proof of 
Theorem 5.
\vskip.3cm
We are now in position to give an expression for the 
M\"obius function $\mu$ of $L(G).$ By Theorem 1 $\mu$ is the 
product of the M\"obius functions of the product lattices.
With the notation introduced thus far, we specifically 
have: 
\vskip.2cm
\noindent
{\bf Theorem 2 \medbreak}
\noindent
{\em (a) The M\"obius function on the lattice of normal }
{\em subgroups of G is} $\mu (H,T)=0,$ {\em unless T/H is in the socle }
{\em S of G/H.}
\medbreak
\noindent
{\em (b) If} $\mbox{\rm $S\cong\prod_{i=1}^aA_i^{d_i}\prod_{j=a+1}^{
a+b}S_j$},$ {\em and }
$T/H\cong\prod_{i=1}^aA_i^{\alpha_i}\prod_{k=1}^cS_{j_k},$ {\em then }
$\mu (H,T)=(-1)^c\prod_{i=1}^a(-1)^{\alpha_i}q_i^{{{\alpha_i}\choose 
2}}.$  
\medbreak
We shall use this form of $\mu$ for explicit calculations 
in Sections 3 and 4 of the paper.
\vskip.5cm
\section*{\Large\bf2. Conjugacy classes that generate a group }
\hspace{\parindent}
A natural question that arises when working on the 
lattice of all subgroups of a group is to identify a 
minimal number of elements that generate the group. 
Hall (1935) calls the cardinality of such a set of elements
the {\em Eulerian number\/} of the group.
The analogous question for the lattice of normal 
subgroups is to consider conjugacy classes that generate 
a group. A set of conjugacy classes is said to {\em generate }
group $G$ if the set of elements in these conjugacy 
classes generate $G.$ 
\vskip.2cm 
\noindent
{\bf Theorem 3} {\em Let} $G$ {\em be a finite group and R be its radical.}
\vskip.1cm
\noindent
(a) {\em If G/R is a direct }
{\em product of non-Abelian simple groups, then there exists }
{\em a conjugacy class that generates G. }
\vskip.1cm
\noindent
(b) {\em If} $C_p^d$ {\em is an elementary Abelian factor of highest dimension }
{\em among all elementary Abelian factors }
{\em in the direct product G/R, then G is generated by d }
{\em conjugacy classes; furthermore, d is the smallest number of }
{\em conjugacy classes required to generate G. }
\vskip.2cm
\noindent
{\bf Proof:}  We often write $\bar {G}$ for $G/R,$ and $\bar {x}$ for the image of 
the element $x\in G$ under the natural epimorphism from $G$ 
to $\bar {G}.$ \vskip.1cm
(a) Let $G/R=S_1\times\cdots\times S_n,$ with $S_i$ non-Abelian simple. 
Think of this group as an abstract group and denote it 
by $S.$ Select $1\neq s_i\in S_i$ and consider the conjugacy class in 
$S$ of the element $s=(s_1,\ldots ,s_n);$ denote this class by 
$s^S.$ Let $H=<s^S>.$ Since $1\neq H\unlhd S$ and the $S_i$ are simple, 
we know that $H\cap S_i,$ being normal in $S_i,$ is either 1 or $
S_i.$ 
By Lemma 1 (b) $H\cap S_i$ cannot be 1 for any $i,$ since 
$1\neq s_i$ is in the projection of $H$ onto $S_i.$ We conclude 
that $H=S(=\bar {G}=G/R).$ View now $s$ as a coset of $R$ in
$G;$ that is, write $s=gR$ for some $g\in G.$ The natural 
epimorphism maps the conjugacy class $g^G$ surjectively 
onto $\bar {g}^{\bar {G}}(=s^S);$ indeed, $\bar {g}^{\bar {y}}
=\overline {g^y.}$ It follows that 
$<g^G,R>=G,$ whence by Lemma 3 (a) we conclude 
$<g^G>=G.$
\vskip.2cm
Turning our attention to part (b), write 
$G/R=S\times (\prod_{j=1}^kC_{p_j}^{d_j}),$ where $S$ is a direct product of 
non-Abelian simple groups and $C_{p_j}^{d_j}$ is an elementary 
Abelian group of dimension $d_j;$ here the $p_j$s are distinct 
primes and $k\geq 1$.  Among the $C_{p_j}^{d_j}$ identify one with maximal $
d_j;$ 
denote it by $C_p^d.$ Write $\bar {G}=G/R$ in the form $\bar {G}=
T\times C_p^d$.
By induction we may assume that $T$ is generated by $m$ 
nonidentity conjugacy classes, with $m\leq d.$ Specifically, let
\[<(x_1,1)^{\bar {G}},\ldots ,(x_m,1)^{\bar {G}}:x_i\in T>=T\times 
1.\]
\noindent
The Abelian subgroup $1\times C_p^d$ is generated by $d$ nonidentity 
conjugacy 
classes, which we write as follows
\[<(1,y_1)^{\bar {G}},\ldots ,(1,y_d)^{\bar {G}}:y_j\in C_p^
d>=1\times C_p^d.\]

\noindent
Consider now the subgroup 
\[H=<(x_1,y_1)^{\bar {G}},\ldots ,(x_m,y_m)^{\bar {G}},\ldots 
,(x_m,y_d)^{\bar {G}}>.\]

\noindent
We assert that $H=\bar {G}.$ To see this, observe that by an 
iterative application of Lemma 1, we conclude that 
$L(\bar {G})=L(T)\times L(C_p^d).$ Consequently $H$ is the direct product 
of its projections on $T$ and on $C_p^d,$ that is 
$H=H_T\times H_{C_p^d}.$ But it is evident from the construction of 
$H$ that $H_T=T$ and $H_{C_p^d}=C_p^d.$ This proves that $H=\bar {
G}.$
Write $\bar {z}_1=(x_1,y_1),\ldots ,\bar {z}_m=(x_m,y_m),\ldots 
,\bar {z}_d=(x_m,y_d),$ with 
$z_i$ being coset representatives. Then 
$<z_1^G,\ldots ,z_d^G,R>=G,$ and therefore $<z_1^G,\ldots ,z_
d^G>=G.$ 
\vskip.2cm
Finally, if fewer than $d$ classes generate $G,$ then 
$\bar {G}=T\times C_p^d$ would 
also be generated by fewer than $d$ classes, namely the 
image of the classes that generate $G$ under the natural 
epimorphism. It follows that the direct factor $C_p^d$ of $\bar {
G}$ 
can be generated by fewer that $d$ elements. This is not 
possible, since $C_p^d$ is a vector space of dimension $d.$ 
End of proof. 
\vskip.2cm
We call the minimal number of classes that generates $G$ 
the {\em class generating number\/} of $G.$ The Eulerian number 
of a group is greater than or equal to the class 
generating number; indeed, if a set of elements generates $G,$ 
than the corresponding conjugacy classes of these 
elements would also generate $G.$
\vskip.5cm
Here are two consequences of Theorem 3. Recall that 
the {\em Frattini\/} subgroup $\phi$ is defined as the intersection of 
all maximal subgroups of a group. In a nilpotent group 
any maximal subgroup is also normal, and this property 
characterizes nilpotent groups. It follows 
that in a nilpotent group $G$ we have $R=\phi ,$ and 
since every maximal
subgroup is of prime index, we conclude that $G/\phi$ is a 
direct product of elementary Abelian groups.
\vskip.2cm
\noindent
{\bf Corollary 1} {\em If G is nilpotent and} $G/\phi =\prod C_{p_
j}^{d_j},$ {\em with} $p_j$ 
{\em distinct primes, then G is generated by} $max_jd_j$ 
{\em elements, and no lesser number of elements generates G.}
\vskip.2cm
\noindent
Indeed, since every maximal subgroup of $G$ is normal, $R$
coincides in this case with the Frattini subgroup of $G.$ 
By Theorem 3 we can generate $G$ with $d=max_jd_j$ 
conjugacy classes. Denote 
these classes by $\{x_i^G:1\leq i\leq d\}.$ If $<x_i:1\leq i
\leq d>\neq G,$ then  
$<x_i:1\leq i\leq d>\leq M,$ for some maximal subgroup $M.$ But 
since $M$ is normal we also have $G=<x_i^G:1\leq i\leq d>\leq 
M,$ 
a contradiction. Minimality of $d$ is also a direct 
consequence of Theorem 3. 
\smallskip\noindent
We recover, in particular, the classical result of 
Burnside, 1899.
\vskip.3cm
\noindent
{\bf Corollary 2} {\em If G is a p-group with a Frattini subgroup }
{\em of index} $p^d,$ {\em then G is generated by d elements, and d }
{\em i\/}s {\em minimal with this property.}
\vskip.1cm
\noindent
Indeed, a $p-$group is nilpotent and in this case $G/\phi =C_p^d.$ The 
statement now follows from Corollary 1.
\vskip.3cm
The class generating number is determined by the 
numerical entries of the character table of the group. On 
the other hand, knowleadge of the class generating 
number places some restrictions on the character table. 
We refer to Fulton and Harris (1991, Part I) for a nicely 
paced memory jogg on group representations. 
To verbalize the interplay, 
always write the character table 
$C$ such that the first row corresponds to the principal 
character, and the first column corresponds to the 
identity conjugacy class. A set of $1+s$ columns of $C$ that 
includes column 1 is called an $s-${\em vertical cut}.  
\vskip.2cm
\noindent
{\bf Proposition 1} {\em A group G has class generating number d }
{\em if and only if d is minimal with the property that }
{\em the character table of G contains a }
{\em d-vertical cut with no row, except the first, having all }
{\em entries equal. }
\vskip.2cm
\noindent
{\bf Proof:} Let $c_1,\ldots ,c_d$ be $d$ (nonidentity) classes that 
generate $G.$ In the 
character table $C$ of $G$ consider the $d-$vertical cut 
corresponding to columns $1,c_1,\ldots ,c_d.$ If row $i$ of this cut 
has equal entries, then the kernel of the irreducible 
representation associated to row $i$ of $C$ contains classes 
$1,c_1,\ldots ,c_d$ and is, therefore, equal to $G;$ this implies 
$i=1.$ Furthermore, assume that there exists a $k-$vertical 
cut having the stated property and that $k<d.$ It follows 
that the $k$ nonidentity classes corresponding to the 
columns of the cut are only in the kernel of the trivial 
representation (associated to the principal character). 
And since every normal subgroup is an intersection of 
kernels of irreducible representations, we conclude that 
the $k$ classes generate the kernel of the trivial 
representation, which is equal to $G;$ this contradicts the 
fact that $d$ is the class generating number of $G.$ 
\vskip.1cm
Conversely, existence of a $d-$vertical cut with minimal 
$d,$ along with the fact that every normal subgroup is an 
intersection of kernels of irreducible representations, 
allows us to conclude that the $d$ classes in the vertical 
cut generate the group $G.$ This ends the proof. 
\vskip.5cm
\section*{\Large\bf3. Class restrictions on the major subgroups}
\hspace{\parindent}
Hopefully in tolerable dissonance with standard terminology, 
we call a normal subgroup of $G$ {\em major\/} if it contains $
R.$
By $\|H\|$ we denote the number of conjugacy classes of 
$G$ contained in the normal subgroup $H.$ The M\"obius 
function used in conjunction with Theorem 3 offers some 
combinatorial insight into the number of conjugacy 
classes contained in the major subgroups of any group.
\vskip.3cm
For $H\unlhd G$ denote by $f_k(H)$ the number of $k-$tuples of 
conjugacy classes of $G$ that generate $H;$ to wit, 
\[f_k(H)=\{(c_1,\ldots ,c_k):<c_1,\ldots ,c_k>=H\}.\]
\noindent
Denote the falling factorial $x(x-1)\cdots (x-m+1)$ by $[x]_
m.$
Since on the lattice $L(G)$ we have $\sum_{T\leq H}f_k(T)=[\|
H\|]_k$, it 
follows by the M\"obius inversion formula that 
\[f_k(G)=\sum_T[\|T\|]_k\mu (T,G),\]
\noindent
with the sum ranging over all major subgroups $T$ of $G.$ 
[This is indeed true, since $\mu (T,G)=\mu (1,G/T)=0,$ unless 
$G/T$ is a direct product of simple groups (by Lemma 2), 
or, equivalently, if $T$ is a major subgroup of $G$ (by 
Lemma 3 (d)).] By the definition of the class generating 
number it follows that $f_k(G)=0,$ for all $k$ less than 
the class generating number of $G.$ We summarize these 
observations below.
\vskip.2cm
\noindent
{\bf Theorem 4 }
\vskip.1cm
\noindent
(a) {\em The number of k-tuples of conjugacy classes that }
{\em generate a group G is equal to}

\[\sum_T[\|T\|]_k\mu (T,G).\]
\vskip.1cm
\noindent
(b) {\em For all} $k$ {\em less than the class generating number of G }
{\em we have}

\[\sum_T[\|T\|]_k\mu (T,G)=0.\]
\noindent
{\em The sum\/}s {\em are over all major subgroups T of G.}
\vskip.2cm
Since Theorem 2 offers explicit values for $\mu ,$ Theorem 
4 (b) provides equational restrictions on the number of 
conjugacy classes contained in the major subgroups of a 
group. 
\vskip.2cm
Examining the case of an Abelian $p-$group, for example,
Theorem 4 specializes to a classical result of Delsarte 
(1948). Indeed, if the type of the Abelian $p-$group $G$ is 
$n_1\geq\cdots\geq n_d,$ with $\sum_in_i=n,$ then $R(G)$ has type
$n_1-1\geq\cdots\geq n_d-1,$ and $G/R(G)$ is a vector space of 
dimension $d$ over the field with $p$ elements. 
Recall that the number of subspaces of 
dimension $i$ in such a vector space is the 
Gaussian polynomial $\big[{d\atop i}\big](p);$ this counts the 
number of major subgroups of order $p^{n-d+i}$ in $G.$ 
Theorem 4 (a) informs us now that the number of 
$k-$tuples of conjugacy classes (in this case elements) that 
generate $G$ is 
\[\sum_{i=0}^d[p^{n-d+i}]_k(\big[{d\atop i}\big](p))(-1)^{d-i}p^{{{
d-i}\choose 2}}.\]
\noindent
As Theorem 4 (b) intimates, this sum is equal to zero 
for all $k$ that are less than $d,$ since $d$ is the class 
generating number (in this instance the number of 
elements in a basis) of $G.$ For $d=k$ the sum counts the 
number of ordered bases of $G.$ Delsarte (1948) finds the 
number of automorphisms, a closely related 
quantity. Since an explicit expression of the M\"obius 
function was not available to him, an indirect method was  
used to determine certain relevant invariants of the 
group. Details on the nonsingularity of the resulting 
coefficient matrix are found in Constantine and 
Kulkarni (1984). As the seminal works of Hall (1933, 1935) 
teach us, the above formula holds in fact in any 
$p-$group of order $p^n$ with Frattini subgroup of index $p^
d.$
Extensions to nilpotent groups are obtained in a 
standard way by working on the $p-$primary parts.
\vskip.5cm
\section*{\Large\bf4. Faithful irreducible representations}
\hspace{\parindent}
Our aim in this section is to examine the column orthogonality relations 
of the character table of a group with respect to only 
the characters of the faithful representations of the group. 
Denote by $\chi_i$ the character of the irreducible representation 
$\rho_i$ of $G,$ over a field $F$ whose characteristic does not 
divide $\vert G\vert .$ Fix conjugacy classes $C$ and $D$ of $G.$ Let 
$f(C,D,N)=\sum^{*}\chi (C)\bar{\chi }(D),$ where the sum is over all irreducible 
characters of $G$ whose 
kernels (of the corresponding $\rho$s) are equal 
to $N$ exactly; the asterix on top of the summation sign 
remins us of that. 
We aim to find an explicit form for $f$ in 
terms of some invariants of $G.$ The {\em indicator function }
$\delta_{AB}(H)$ is by definition equal to 0 if conjugacy classes $
A$ 
and $B$ of group $H$ are distinct in $H$, and is equal to 1 if 
they are the same in $H.$ To simplify notation we write 
$\delta_{AB}(G/K)$ for $\delta_{\bar {A}\bar {B}}(G/K)$, where $A$ and $
B$ are classes of $G$ 
and $\bar {A}$ and $\bar {B}$ are their epimorphic images in the factor 
group $G/K.$
\medbreak
Using the notation from Theorems 1 and 2, we prove 
the following result.
\medbreak
\noindent
{\bf Theorem 5 }
\medbreak
\noindent
{\em (a) If conjugacy classes A and B are distinct modulo }
{\em the socle S of group G, then }
\[\sum^{*}\chi (A)\bar{\chi }(B)=0.\]
\medbreak
\noindent
{\em (b) If conjugacy class C is such that} $\vert CS\vert =\vert 
C\vert\vert S\vert ,$ {\em then}
\[\sum^{*}\chi (C)\bar{\chi }(C)=\frac {\vert G\vert}{\vert C\vert}
\prod_{i=1}^a\prod_{j=0}^{d_i-1}(1-\frac {q_i^j}{\vert A_i\vert})
\prod_{k=a+1}^{a+b}(1-\frac 1{\vert S_k\vert}).\]
{\em (c) The sum of squared dimensions of the faithful }
{\em irreducible representations, i.e.} $\sum^{*}\chi (1)^2${\em , is divisible by }
{\em the cardinality of any conjugacy class of G that }
{\em satisfies} $|CS|=|C||S|.$
\vskip.5cm
\medbreak\noindent
Before proving Theorem 5 we make a few remarks. To simplify the language, call classes $
A$ and $B$ {\em faithfully }
{\em orthogonal\/} if $\sum^{*}\chi (A)\chi (B)=0.$ Part (a) of Theorem 5 than 
states that {\em Any two conjugacy classes that are distinct modulo the }
{\em socle are faithfully orthogonal.\/} In particular, {\em Any }
{\em conjugacy class inside the socle is faithfully }
{\em orthogonal to any conjucagy class outside the socle. }
Furthermore, it also immediately follows that {\em If m }
{\em denotes the number of conjugacy classes of G/S, then  }
{\em there are at least m conjugacy classes of group G any }
{\em (distinct) two of which are faithfully orthogonal.\/} A 
special case of note, P\'alfy (1979), occurs in Part (b) when the conjugacy 
class $C=1;$ the sum $\sum^{*}\chi (C)\bar{\chi }(C)$ is in such a case equal to 
the sum of squares of the dimensions of the faithful 
irreducible representations. We observe that the condition 
$|CS|=|C||S|$ is easily seen to always be satisfied when $C$ is in the 
center of $G,$ but many noncentral such classes exist, in 
general. [For instance, in the group $SL(2,3)$. There is a single minimal 
normal subgroup, namely $S=<\left[\begin{array}{cc}
-1&0\\
0&-1\end{array}
\right]>$, which is the socle
of $SL(2,3).$ Theorem 5 (b) yields 
$\sum_f\chi (1)^2=24(1-\frac 12)=12,$ since in this case $a=1,$ $
b=0,$ 
and $d_1=1.$ Consider the conjugacy class $C$ of the element 
$\left[\begin{array}{cc}
1&1\\
0&1\end{array}
\right]$. Class $C$ contains 4 elements, and it satisfies 
$\vert CS\vert =\vert C\vert\vert S\vert .$ Indeed, $SL(2,3)$ has three faithful 
irreducible representations of degree 2, and 
$12=\sum_f\chi (1)^2=2^2+2^2+2^2$ is divisible by $4=\vert C\vert 
,$ as 
Theorem 5 (c) states.] In fact, for certain groups, 
the divisibility stated in Theorem 5 (c) may prove 
useful in sheding light on the invariants $a,$$b,$ and $d_i$ of 
the socle, which appear in Theorem 5 (b).
\medbreak
We proceed with the proof of Theorem 5. By the column 
orthogonality of the character table of $G,$ we have
\[\sum_{N\leq T}f(A,B,T)=(|G|/|AN|)\delta_{AB}(G/N),\]
where the sum is over all irreducible representations of 
$G$ whose kernel contain $N.$ Interpreting this equality on 
the lattice $L(G)$, we use M\"obius inversion to obtain
\[f(A,B,N)=\sum_{N\leq T}(|G|/|AT|)\delta_{AB}(G/T)\mu (N,T).\]
By working in the group $G/N$ it suffices to evaluate 
$f(A,B,1)$. Using Hall's result we may further simplify by 
restricting the sum only over normal subgroups $T$ of $G$ 
that are in the socle $S$ of $G.$ These observations yield 
\[f(A,B,1)=\sum_{T\leq S}(|G|/|AT|)\delta_{AB}(G/T)\mu (1,T).\]
The last expression simplifies to closed form in certain 
instances. Indeed, it is easy to see that if 
$\delta_{AB}(G/S)=0,$ then $\delta_{AB}(G/T)=0$ for all $T\leq S.$ We 
therefore conclude that when classes $A$ and $B$ are 
distinct modulo the socle $S,$ then every term in this 
expression is, in fact, zero. This proves Part (a). 
\smallskip\noindent
Turning attention to Part (b), in case $A=B$ to simplify notation write 
$f(A,N)$ for $f(A,B,N).$ Let $C$ be a conjugacy class of $G.$ If $
N$ is a normal subgroup of $G,$ then by the column 
orthogonality of the character table of $G,$ we have 
\[\sum_{N\leq T}f(C,T)=\vert G\vert /\vert CN\vert ,\]
where the sum is over all irreducible representations 
whose kernels contain $N.$ Iinvert to 
obtain
\[f(C,N)=\sum_{N\leq T}(\vert G\vert /\vert CT\vert )\mu (N,T).\]
Since in Theorem 2 we obtained an expression for $\mu ,$ 
the right hand side of this equality can now be, in 
principle at least, 
explicitely computed. By working, again, in the group $G/N$ it 
suffices to evaluate $f(C,1).$ Following the notation of Theorem 
2 (with $H=1$) we obtain:
\medbreak
 
\[f(C,1)=\sum_{T\leq S}(\vert G\vert /\vert CT\vert )\mu (1,T).\]
\smallskip\noindent
In certain instances this expression can be written as 
a product. Assume that the class $C$ is such that 
$\vert CS\vert =\vert C\vert\vert S\vert ,$ where $S$ is the socle of $
G;$ that is, class $C$ 
preserves its cardinality when passing to the group $G/S.$ 
For any $T\leq S,$ $T$ normal in $G,$ we 
then obviously have $\vert CT\vert =\vert C\vert\vert T\vert ,$ since each coset of $
S$ 
contains at most one element of $C$, and therefore the 
same must be true of cosets of $T.$ We thus obtain,
\[f(C,1)=\sum_{T\leq S}(\vert G\vert /\vert C\vert\vert T\vert )\mu 
(1,T)\]
\[=\frac {\vert G\vert}{\vert C\vert}\sum_{T_1\leq S_1}\cdots\sum_{
T_{a+s}\leq S_{a+s}}(1/\prod T_i\vert )\mu (1,\prod T_i)\]
\[=\frac {\vert G\vert}{\vert C\vert}\prod_{i=1}^{a+s}(\sum_{T_i\leq 
S_i}(1/\vert T_i\vert )\mu (1,T_i)).\]
\smallskip
\noindent
Using Theorems 1 and 2, we can turn the above sums 
over the elements of a lattice into ordinary sums. For 
$a+1\leq i\leq a+b$, we have $\sum_{T_i\leq S_i}(1/\vert T_i\vert 
)\mu (1,T_i)=1-\frac 1{\vert S_i\vert}$, 
since in this case $S_i$ is nonabelian minimal normal in $G$, 
and therefore $T_i$ is either 1 or $S_i.$ As to the Abelian 
case, for $1\leq i\leq a$ we have
\[\sum_{T_i\leq S_i}\vert T_i\vert^{-1}\mu (1,T_i)=\sum_{k=0}^{d_
i}\vert A_i^{}\vert^{-k}\frac {(q_i^{d_i}-1)\cdots (q_i^{d_i}-q_i^{
k-1})}{(q_i^k-1)\cdots (q_i^k-q_i^{k-1})}(-1)^kq_i^{{k\choose 2}}\]
\[=\sum_{k=0}^{d_i}\vert A_i^{}\vert^{-k}(-1)^k\frac {(q_i^{d_i}-
1)\cdots (q_i^{d_i}-q_i^{k-1})}{(q_i^k-1)\cdots (q_i-1)}.\]
\smallskip
\noindent
Using the $q-$identity
\[\frac {(q^m-1)\cdots (q^m-q^{k-1})}{(q^k-1)\cdots (q-1)}+q_{}^m\frac {
(q^m-1)\cdots (q^m-q^{k-2})}{(q^{k-1}-1)\cdots (q-1)}\]
\[=\frac {(q^m-1)\cdots (q^m-q^{k-2})}{(q^{k-1}-1)\cdots (q-1)}(\frac {
q^m-q^{k-1}}{q^k-1}+\frac {q^{m+k}-q^m}{q^k-1})\]
\[=\frac {(q^{m+1}-1)\cdots (q^{m+1}-q^{k-1})}{(q^k-1)\cdots (q-1
)},\]
\smallskip
\noindent
we obtain
\[(\sum_{k=0}^m\vert A_i^{}\vert^{-k}(-1)^k\frac {(q_i^m-1)\cdots 
(q_i^m-q_i^{k-1})}{(q_i^k-1)\cdots (q_i-1)})(1-\frac {q_i^m}{\vert 
A_i\vert})\]
\[\sum_{k=0}^{m+1}\vert A_i^{}\vert^{-k}(-1)^k\frac {(q_i^{m+1}-1
)\cdots (q_i^{m+1}-q_i^{k-1})}{(q_i^k-1)\cdots (q_i-1)},\]
\smallskip
\noindent
and hence, by induction,
\[\sum_{k=0}^{d_i}\vert A_i^{}\vert^{-k}(-1)^k\frac {(q_i^{d_i}-1
)\cdots (q_i^{d_i}-q_i^{k-1})}{(q_i^k-1)\cdots (q_i-1)}=\prod_{j=
0}^{d_i-1}(1-\frac {q_i^j}{\vert A_i\vert}).\]
\smallskip
\noindent
Therefore,
\[f(C,1)=\frac {\vert G\vert}{\vert C\vert}\prod_{i=1}^a\prod_{j=
0}^{d_i-1}(1-\frac {q_i^j}{\vert A_i\vert})\prod_{k=a+1}^{a+b}(1-\frac 
1{\vert S_k\vert}),\]
which proves Part (b).
\medbreak\noindent
As to Part (c), observe that Theorem 5 (b) can be writen as
$\sum^{*}\chi (C)\bar{\chi }(C)=\frac 1{\vert C\vert}\sum^{*}\chi 
(1)^2.$ The right hand side is  
a rational number, whereas the left hand side 
is an algebraic integer. It follows that the right hand 
side is in fact a rational integer (possibly zero). This ends the 
proof of Theorem 5.
\vskip.3cm
\noindent
An immediate consequence of Theorem 5 is the 
following:
\vskip.2cm
\noindent
{\bf Corollary} 3 {\em Group G does not have a faithful }
{\em irreducible representation if and only if for some i, }
$1\leq i\leq a,$ {\em and some j,} $0\leq j\leq d_i-1${\em ,} $\vert 
A_i\vert =\vert Hom_G(A_i,A_i)\vert^j.$
\vskip.3cm
\noindent
We examine Corollary 3 when the group $G$ is nilpotent.
Then $G=P_1\times\cdots\times P_n$, where $P_i$ are the Sylow 
subgroups of $G.$ If $S(H)$ denotes the socle of group $H,$ 
then it is easy to verify that $S(G)=S(P_1)\times\cdots\times S(P_
n)$.
Consider a $p-$group $P.$ The class equation informs us that 
any minimal normal subgroup of $P$ intersects the center $Z(P)$ of $
P$ 
nontrivially, and is therefore entirely contained in $Z(P).$ 
It follows that $S(P)\leq Z(P).$ There is, therefore, only one 
$G-$isomorphism type of minimal normal subgroups of $P,$ 
namely a subgroup $C_p$ of order $p$ which $G$ fixes pointwise. 
In this case $\vert Hom_P(C_p,C_p)\vert =p,$ and $S(P)\cong C_p^d$ for some 
$d\geq 1.$ The above Corollary informs us that $P$ has a 
faithful irreducible representation if and only if 
$d=1,$ and this happens (by the basis theorem for Abelian 
$p-$groups) if and only if $Z(P)$ is cyclic. Using commutation 
of
$Hom$ with direct products, and working on the $p-$primary 
parts $P_i$ of $G,$ we conclude that $G$ has a faithful 
irreducible representation if and only if $d_i=1$, for all i. 
The last condition is equivalent to saying that 
$Z(G)=Z(P_1)\times\cdots\times Z(P_n)$ is a cyclic group. We thus 
argued from a different perspective the validity of 
the following well-known result; see Gorenstein (1968, p 64):
\vskip.2cm
\noindent
{\bf Corollary} 4 {\em A nilpotent group has a faithful irreducible }
{\em representation if and only if it has a cyclic center.}
\vfill\eject
\centerline{REFERENCES}
\begin{enumerate}

\item Alperin, J. L. and Bell, R. B. (1995) {\em Groups and }
{\em representations,\/} Springer, New York.

\item Burnside, W. (1911) {\em Theory of groups,} $2^{nd}$ edition, 
Dover (New York, 1955).

\item Constantine, G. (1987).  {\it Combinatorial Theory and
Statistical Design,} Wiley, New York. 

\item Constantine, G. and Kulkarni, R. S. (1984) On a result 
of S. Delsarte, {\em Proc. AMS,\/} {\bf 92}, 149-152.

\item Delsarte, S. (1948) Fonctions de M\"obius sur les 
groupes Abeliens finis, {\em Ann. Math.,\/} {\bf 49}, 600-609.

\item Fulton, W. and Harris, J. (1991) {\em Representation }
{\em Theory,\/} Springer, New York.

\item Gorenstein, D. (1968) {\em Finite Groups}, Harper and Row, 
New York.

\item Hall, P. (1933) A contribution to the theory of 
groups of prime-power order, {\em Proc. London Math. Soc., }
{\bf 36}, 29-95. 

\item Hall, P. (1935) The Eulerian function of a group, 
{\em Quart. J. Math.,\/} {\bf 7}, 134-151.

\item Miller, M. D. (1975) On the lattice of normal subgroups of a
direct product, {\em Pacific J. Math.}, 60, no. 2, 153-158.

\item P\'alfy, P. (1979) On faithful irreducible 
representations of finite groups, {\em Studia Sci. Math. }
{\em Hungar.,\/} {\bf 14}, 95-98.

\item Rota, G-C. (1964) On the foundations of combinatorial 
theory, I. Theory of M\"obius functions, {\em Z. Wahrsch. Verw. }
{\em Gebiete,\/} {\bf 2}, 340-368.

\item Silcock, H. L. (1977) Generalized wreath products and the
lattice of normal subgroups of a group,
{\em Algebra-Universalis}, 3, 361-372.

\item Stanley, R. (1982) Some aspects of groups acting on finite 
posets, {\em J. Combinatorial Theory (A)}, {\bf 32}, (1982), 132-161.

\item Weisner, L. (1935) Abstract theory of inversion of 
finite series, {\em Trans. Amer. Math. Soc.,\/} {\bf 38}, 474-484.
\end{enumerate}
\end{document}